\documentclass[a4paper, 9pt, reqno]{amsart}
\usepackage[varg]{txfonts}
\usepackage{verbatim}
\usepackage{amsfonts, amsmath, amssymb, amsthm, enumerate, color}
\usepackage[left=1in, right=1in, top=0.9in, bottom=0.8in]{geometry}

\newcommand{\R}{\mathbb{R}}

\newtheorem{thm}{Theorem}[section]

\newtheorem{prop}{Proposition}[section]

\author[J. Seok]{Jinmyoung Seok}
\address[Jinmyoung Seok]{Department of Mathematics, Kyonggi University,
154-42 Gwanggyosan-ro, Yeongtong-gu, Suwon 16227, Republic of Korea}
\email{jmseok@kgu.ac.kr}

\title[Nonlinear Choquard equations: doubly critical case]
{Nonlinear Choquard equations: doubly critical case}

\begin{document}

\maketitle

\begin{abstract}
Consider nonlinear Choquard equations
\begin{equation*}
\left\{\begin{array}{rcl}
-\Delta u +u & = &(I_\alpha*F(u))F'(u)  \quad \text{in } \R^N, \\
\lim_{x \to \infty}u(x) & = &0,
\end{array}\right.
\end{equation*} 
where $I_\alpha$ denotes Riesz potential and $\alpha \in (0, N)$. 
In this paper, we show that when $F$ is doubly critical, i.e.
$F(u) = \frac{N}{N+\alpha}|u|^{\frac{N+\alpha}{N}}+\frac{N-2}{N+\alpha}|u|^{\frac{N+\alpha}{N-2}}$,
the nonlinear Choquard equation admits a nontrivial solution if $N \geq 5$ and $\alpha + 4 < N$.
\end{abstract}

\section{Introduction}
Let $N \geq 3,\, \alpha \in (0, N)$.
We are concerned with the nonlinear Choquard equation:
\begin{equation}\label{main-eq}
\left\{\begin{array}{rcl}
-\Delta u +u & = & (I_\alpha*F(u))F'(u)
\quad \text{in } \R^N, \\
\lim_{x \to \infty}u(x) & = &0,
\end{array}\right.
\end{equation}
where $I_\alpha$ is Riesz potential given by
\[
I_\alpha(x) = \frac{\Gamma(\frac{N-\alpha}{2})}{\Gamma(\frac{\alpha}{2})\pi^{N/2}2^\alpha|x|^{N-\alpha}}
\]
and $\Gamma$ denotes the Gamma function.
It is the Euler-Lagrange equation of the functional
\[
J_\alpha(u) = \frac12\int_{\R^N}|\nabla u|^2+u^2 \,dx -\frac{1}{2}\int_{\R^N}(I_\alpha*F(u))F(u)\,dx.
\]
Physical motivation of \eqref{main-eq} comes from the case that $F(u) = \frac{1}{2}|u|^2$ and $\alpha = 2$.
In this case, the equation \eqref{main-eq} is called the Choquard-Pekar equation \cite{L, P},
Hartree equation \cite{FL, LNR} or Schr\"odinger-Newton equation \cite{MPT, TM}, depending on its physical backgrounds and derivations.  
The existence of a ground state in this case is studied in \cite{L, Li, M} via variational arguments.

The functional $J_\alpha$ can be considered as a nonlocal perturbation of the fairly well-studied functional consisting of only local terms:
\[
J_0(u) = \frac12\int_{\R^N}|\nabla u|^2+u^2 \,dx -\int_{\R^N}G(u)\,dx
\] 
since as $\alpha \to 0$, $J_\alpha$ approaches to $J_0$ with $G(u) = \frac12F^2(u)$.
A critical point of $J_0$ is a solution to the stationary nonlinear Schr\"odinger equation:
\begin{equation}\label{limit-eq}
-\Delta u + u = G'(u). 
\end{equation}

The power type function $\frac{1}{p}|u|^p$ is a standard choice for nonlinearity $G(u)$ (and also $F(u)$). 
By Sobolev inequality, it can be shown that the functional $J_0$ is a well-defined $C^1$ functional on $H^1(\R^N)$ if $G(u) = \frac{1}{p}|u|^p$ and $p \in [2, \frac{2N}{N-2}]$. 
It is a classical result that it admits a nontrivial critical point of ground state level in the subcritical range $p \in (2, \frac{2N}{N-2})$ \cite{BL, Strauss}. 
Moreover, the standard application of Pohozaev's identity says that if $p$ is out of subcritical, i.e., $1 < p \leq 2$ or $p \geq 2N/(N-2)$, the equation \eqref{limit-eq} does not admit any nontrivial finite energy solution. 
In case of $J_\alpha$, Hardy-Littlewood-Sobolev inequality (Proposition \ref{HLS} below) replaces Sobolev inequality to see that 
$J_\alpha$ with $F(u) = \frac{1}{p}|u|^p$ is well-defined and is continuously differentiable on $H^1(\R^N)$ if $p \in [\frac{N+\alpha}{N},\, \frac{N+\alpha}{N-2}]$. 
Two numbers $\frac{N+\alpha}{N}$ and $\frac{N+\alpha}{N-2}$ play roles of lower and upper critical exponents for existence. 
It is proved by Moroz and Van schaftingen \cite{MV2} that for every $\alpha \in (0,\, N)$, there exists a nontrivial ground state solution  
if $p$ is in the subcritical range, i.e., $p \in (\frac{N+\alpha}{N},\, \frac{N+\alpha}{N-2})$
and there is no nontrivial finite energy solution  if $p$ is outside of subcritical, i.e.,  $1 < p \leq \frac{N+\alpha}{N}$ or $p \geq \frac{N+\alpha}{N-2}$.
This result is compatible with the existence of the limit equation \eqref{limit-eq}. 
Observe the existence range $p \in (\frac{N+\alpha}{N},\, \frac{N+\alpha}{N-2})$ tends to $p \in (1,\, \frac{N}{N-2})$ and the nonlinear term $(I_\alpha*|u|^p)|u|^p$ tends to $|u|^{2p}$   as $\alpha \to 0$. 
We recall that the equation \eqref{limit-eq} with $G(u) = \frac{1}{2}F^2(u) = \frac{1}{2p^2}|u|^{2p}$ admits a nontrivial finite energy solution if and only if $p \in (1, \frac{N}{N-2})$.
Furthermore, we have $H^1$ convergence between ground states. For any $p \in (1, \, \frac{N}{N-2})$, choose a small $\alpha_0 > 0$ 
that $p$ belongs to the segment $(\frac{N+\alpha}{N},\,\frac{N+\alpha}{N-2})$ for every $\alpha \in (0,\, \alpha_0)$ so that a radial positive ground state $u_\alpha$ to \eqref{main-eq} with $F(u) = \frac{1}{p}|u|^{p}$ exists. 
Then it is possible to show that as $\alpha \to 0$, $u_\alpha$ converges in $H^1$ sense to a ground state $u_0$ of the corresponding functional  $J_0$. See \cite{RV, S}.

For general nonlinearity $G$, Berestycki and Lions prove in their celebrated paper \cite{BL} that \eqref{limit-eq} admits a ground state solution when $G$ is $C^1(\R)$ and satisfies the following:
\begin{enumerate}[(G1)]
\item there exists a constant $C > 0$ such that for every $s \in \R$, 
\[
|sG'(s)| \leq C(|s|^2 +|s|^{\frac{2N}{N-2}}),
\] 
\item 
$\displaystyle{\lim_{s\to\infty}\frac{G(s)}{|s|^\frac{2N}{N-2}} = 0}\quad$ and 
$\quad\displaystyle{\lim_{s\to0}\frac{G(s)}{|s|^2} = 0}$, 
\\
\item there exists a constant $s_0 \in \R \setminus \{0\}$ such that $G(s_0) > \frac{s_0^2}{2}$. \smallskip
\end{enumerate}

In the same spirit, it is proved in \cite{MV1} that there exists a ground state solution to \eqref{main-eq} under the following conditions for the nonlinearity function $F \in C^1(\R)$: \smallskip
\begin{enumerate}[(F1)]
\item (growth) there exists a constant $C > 0$ such that for every $s \in \R$, 
\[
|sF'(s)| \leq C(|s|^{\frac{N+\alpha}{N}} +|s|^{\frac{N+\alpha}{N-2}}), 
\] 
\item (subcriticality)
$\displaystyle{\lim_{s\to\infty}\frac{F(s)}{|s|^\frac{N+\alpha}{N-2}} = 0}\quad$ and $\quad\displaystyle{\lim_{s\to0}\frac{F(s)}{|s|^{\frac{N+\alpha}{N}}} = 0}$,
\\
\item (nontriviality) there exists a constant $s_0 \in \R \setminus \{0\}$ such that $F(s_0) \neq 0$. \smallskip
\end{enumerate}
\noindent Interestingly, note that the condition (F3) are inconsistent with (G3) of limit equation \eqref{limit-eq}
while we have seen the consistency between \eqref{main-eq} and \eqref{limit-eq} of power type.

In this paper, we are interested in some choice of $F$ that violates the subcriticality condition (F2). The setting that we can naturally choose would be
\[
F(u) \coloneqq \frac{1}{p}|u|^{p} +\frac{1}{q}|u|^{q}, \quad \frac{N+\alpha}{N} \leq p <  q \leq \frac{N+\alpha}{N-2}
\]
and either $p = \frac{N+\alpha}{N}$ or $q =  \frac{N+\alpha}{N-2}$.
With $G(u) = \frac1p|u|^p+\frac1q|u|^q$, it is proved in \cite{ADS} that the limit equation \eqref{limit-eq} admits a nontrivial solution if $2 < p < q = \frac{2N}{N-2}$ for $N \geq 4$ and $4 < p < 6$ for $N = 3$.
In this setting, considering a doubly critical choice of $p$ and $q$, i.e., $p = 2,\, q = \frac{2N}{N-2}$ does not seem appropriate because the equation is just reduced to the Lane-Emden equation $-\Delta u = |u|^{\frac{4}{N-2}}u$.
The situation is however different when we consider the nonlinear Choquard equation \eqref{main-eq} with a pair of lower and upper critical exponents: $p = \frac{N+\alpha}{N}$ and $q = \frac{N+\alpha}{N-2}$. 
The purpose of this paper is to study this case.  
We shall prove that there exists a nontrivial solution under some restrictions on $N$ and $\alpha$.
\begin{thm}\label{doubly-critical}
Let $N \geq 5$ and $F(u) = \frac1p|u|^p +\frac1q|u|^q$.
Then, there exists a nontrivial solution $u \in H^1_r(\R^N)$ to \eqref{main-eq} if 
$p = \frac{N+\alpha}{N},\, q = \frac{N+\alpha}{N-2}$ and $N > 4+\alpha$.
\end{thm}

For the related critical problems involving only a single critical exponent, we refer to \cite{AGSY, Ao, CZ, MV3, S2}.
When we approach by variational methods to prove Theorem \ref{doubly-critical}, the main difficulty we encounter is to deal with two different types of loss of compactness. 
By expanding the nonlinear term $\int_{\R^N}(I_\alpha*F(u))F(u)$, we see that it contains two terms
$\int_{\R^N}(I_\alpha*|u|^\frac{N+\alpha}{N})|u|^\frac{N+\alpha}{N}$ and $\int_{\R^N}(I_\alpha*|u|^\frac{N+\alpha}{N-2})|u|^\frac{N+\alpha}{N-2}$, which are invariant under dilations $\lambda^{N/2}u(\lambda\cdot)$ and $\lambda^{(N-2)/2}u(\lambda\cdot)$ respectively. 
These two dilations are noncompact group actions on $H^1(\R^N)$, each of which prevents a general (PS) sequence of $J_\alpha$ from being relatively compact.
The following two inequality, that is special cases of Hardy-Littlewood-Sobolev inequality, play significant roles to resolve this difficulty.
The first one is
\begin{equation}\label{HLS1}
\mathcal{S}_1\left(\int_{\R^N}(I_\alpha*|u|^{\frac{N+\alpha}{N}})|u|^{\frac{N+\alpha}{N}}\right)^{\frac{N}{N+\alpha}}
\leq \int_{\R^N}u^2\,dx
\end{equation}
whose extremal functions are 
\[
u(x) = C\frac{\lambda^{N/2}}{(\lambda^2+|x|^2)^{N/2}}.
\]
The second one is
\begin{equation}\label{HLS2}
\mathcal{S}_2\left(\int_{\R^N}(I_\alpha*|u|^{\frac{N+\alpha}{N-2}})|u|^{\frac{N+\alpha}{N-2}}\right)^{\frac{N-2}{N+\alpha}}
\leq \int_{\R^N}|\nabla u|^2\,dx
\end{equation}
whose extremal functions are
\[
u(x) = C\frac{\lambda^{(N-2)/2}}{(\lambda^2+|x|^2)^{(N-2)/2}}.
\]

A key step to prove the existence of a solution is a characterization of level sets at which a (PS) sequence of $J_\alpha$ converges.
These levels are given in terms of two best constants $\mathcal{S}_1$ and $\mathcal{S}_2$ of inequalities \eqref{HLS1} and \eqref{HLS2}.
More precisely, we shall obtain the following proposition. We say a sequence $\{u_j\} \subset H^1_r(\R^N)$ is a (PS) sequence of $J_\alpha|_{H^1_r(\R^N)}$ at level $c$ if
\[
J_\alpha'(u_j) \to 0 \text{ in } H^{-1}_r(\R^N) \quad \text{and} \quad J_\alpha(u_j) \to c \quad \text{as } j \to \infty.
\]
\begin{prop}\label{compact-level}
Assume $p = \frac{N+\alpha}{N},\, q = \frac{N+\alpha}{N-2}$.
Let $\{u_j\} \subset H^1_r(\R^N)$ be a $(PS)$ sequence of $J_\alpha|_{H^1_r(\R^N)}$ at level $c$.
Then it is relatively compact in $H^1(\R^N)$ if 
\[
c < \min\left(\frac{1}{2}(1-\frac{1}{p})p^{\frac{1}{p-1}}\mathcal{S}_1^{\frac{p}{p-1}},\,   
\frac{1}{2}(1-\frac{1}{q})q^{\frac{1}{q-1}}\mathcal{S}_2^{\frac{q}{q-1}}\right).
\]
\end{prop}
Proposition \ref{compact-level} shall be proved in Section 3. 
Then, the remaining work is to show the mountain pass energy level $c$ of $J_\alpha$ satisfies the condition of Proposition \ref{compact-level}.  
This shall be done in Section 4 by testing two aforementioned families of extremal functions to the functional $J_\alpha$.
In Section 2, we collect various useful inequalities and estimate required when we prove Proposition \ref{compact-level} and Theorem \ref{doubly-critical}.

%


\section{Auxiliary tools}
In this section, we prepare some auxiliary tools for proving our main theorem.
The well-known Hardy-Littlewood-Sobolev inequality is stated as follows.
\begin{prop}[Hardy-Littlewood-Sobolev inequality \cite{Gra, LL}]\label{HLS}
Let $p, r > 1$ and $0 < \alpha < N$ be such that
\[
\frac{1}{p}+\frac{1}{r} = 1+\frac{\alpha}{N}.
\]
Then there exists $C > 0$ depending only on $N,\alpha,p$ such that for any $f \in L^p(\R^N)$ and $g \in L^r(\R^N)$
\[
\left|\int_{\R^N}\!\!\int_{\R^N}\frac{f(x)g(y)}{|x-y|^{N-\alpha}}\,dxdy\right| \leq C(N,\alpha,p)\|f\|_{L^p(\R^N)}\|g\|_{L^r(\R^N)}.
\]
%
\end{prop}

Hardy-Littlewood-Sobolev inequality has a dual form called the Riesz potential estimate.
\begin{prop}[Riesz potential estimate \cite{Gra, LL}]\label{est-riesz}
Let $1 \leq p < s < \infty$ and $0 < \alpha < N$ be such that
\[
\frac1p-\frac1s = \frac{\alpha}{N}.
\]
Then there exists $C > 0$ depending only on $N,\alpha,p$ such that for any $f \in L^p(\R^N)$,
\[
\left\|\frac{1}{|\cdot|^{N-\alpha}}*f\right\|_{L^s(\R^N)} \leq C \|f\|_{L^p(\R^N)}.
\]
\end{prop}

We denote by $H^1_r(\R^N)$ the space of radial functions in $H^1(\R^N)$.
Non-invariance of $H^1$ norm of a function $u \in H^1_r(\R^N)$ by translations induces the compact embedding to $L^p(\R^N)$ for subcritical $p$. See \cite{Strauss}.
\begin{prop}
The Sobolev embedding $H^1_r(\R^N) \hookrightarrow L^p(\R^N)$ is compact if $2 < p < 2N/(N-2)$.
\end{prop}

By combining Hardy-Littlewood-Sobolev inequality (Proposition \ref{HLS}) and the compact Sobolev embedding, a standard analysis shows the following convergences hold. We refer to \cite{S} for details.
\begin{prop}\label{convergence}
Let $\alpha \in (0,\, N)$ and $\{u_j\} \subset H_r^1(\R^N)$ be a sequence converging weakly to some $u_0 \in H^1_r(\R^N)$ in $H^1(\R^N)$ as $j \to \infty$.
\begin{enumerate}[$(i)$]
\item If $\frac{N+\alpha}{N} < p \leq q < \frac{N+\alpha}{N-2}$, then
\[
\int_{\R^N}\left(\frac{1}{|\cdot|^{N-\alpha}}*|u_j|^p\right)|u_j|^q\,dx 
\to \int_{\R^N}\left(\frac{1}{|\cdot|^{N-\alpha}}*|u_0|^p\right)|u_0|^q\,dx;
\]
\item If $\phi \in H^1(\R^N)$, $\frac{N+\alpha}{N} \leq p \leq q \leq \frac{N+\alpha}{N-2}$, then
\[
\int_{\R^N}\left(\frac{1}{|\cdot|^{N-\alpha}}*|u_j|^p\right)|u_j|^{p-2}u_j\phi\,dx 
\to \int_{\R^N}\left(\frac{1}{|\cdot|^{N-\alpha}}*|u_0|^p\right)|u_0|^{p-2}u_0\phi\,dx.
\]
\end{enumerate}
\end{prop}

The following version of Brezis-Lieb lemma for the Riesz potential is useful for our analysis. We refer to \cite{MV2} for a proof.
\begin{prop}\label{Brezis-Lieb}
Let $\alpha \in (0,\, N)$ and $p \in [\frac{N+\alpha}{N},\, \frac{N+\alpha}{N-2}]$ be given.
If $\{u_n\}$ be a bounded sequence in $L^{\frac{2Np}{N+\alpha}}(\R^N)$ such that
$u_n \to u$ almost everywhere as $n \to \infty$ for some function $u$, then 
$u \in L^{\frac{2Np}{N+\alpha}}(\R^N)$ and
\[
\lim_{n\to\infty}\left(\int_{\R^N}(I_\alpha*|u_n|^p)|u_n|^p-\int_{\R^N}(I_\alpha*|u_n-u|^p)|u_n-u|^p\right)
= \int_{\R^N}(I_\alpha*|u|^p)|u|^p.
\] 
\end{prop}

\section{Proof of Proposition \ref{compact-level}}

In this section, we prove Proposition \ref{compact-level}.
We first show that $\{u_j\}$ is bounded in $H^1(\R^N)$.
Indeed, from the definition of $(PS)$ sequence,
\[
\begin{aligned}
&c+o(1) = J_\alpha(u_j) = \frac{1}{2}\|u_j\|_{H^1}^2-\frac{1}{2}\int_{\R^N}(I_\alpha*(\frac{1}{p}|u_j|^p+\frac{1}{q}|u_j|^q))(\frac{1}{p}|u_j|^p+\frac{1}{q}|u_j|^q)\,dx,  \\
&o(1)\|u_j\|_{H^1} = J'_\alpha(u_j)u_j = \|u_j\|_{H^1}^2-\int_{\R^N}(I_\alpha*(\frac{1}{p}|u_j|^p+\frac{1}{q}|u_j|^q))(|u_j|^p+|u_j|^q)\,dx.
\end{aligned}
\]
Then,

\[
\begin{aligned}
\frac{1}{2}\|u_j\|_{H^1}^2 & \leq c+o(1) +\frac{1}{2p}\int_{\R^N}(I_\alpha*(\frac{1}{p}|u_j|^p+\frac{1}{q}|u_j|^q))(|u_j|^p+|u_j|^q)\,dx \\
& = \frac{1}{2p}\left(\|u_j\|_{H^1}^2 +o(1)\|u_j\|_{H^1}\right) +c +o(1).
\end{aligned}
\]
Since $p > 1$, this shows $\|u_j\|_{H^1}$ is bounded.

Now, up to a subsequence, $\{u_j\}$ weakly converges to some $u_0 \in H^1_r(\R^N)$.
Using $(ii)$ of Proposition \ref{convergence}, it is standard to show that $u_0$ is a weak solution of \eqref{main-eq}.
Let $w_j := u_j-u_0$. 
From Proposition \ref{Brezis-Lieb} and Proposition \ref{convergence}, we see that
\begin{equation}\label{equ2}
\begin{aligned}
\|w_j\|_{H^1}^2 &= \|u_j-u_0\|_{H^1}^2 = \|u_j\|_{H^1}^2-\|u_0\|_{H^1}^2 +o(1) \\
&= \int_{\R^N}(I_\alpha*(\frac{1}{p}|u_j|^p+\frac{1}{q}|u_j|^q))(|u_j|^p+|u_j|^q)\,dx +o(1)\|u_j\|_{H^1} \\
&\qquad\qquad-\int_{\R^N}(I_\alpha*(\frac{1}{p}|u_0|^p+\frac{1}{q}|u_0|^q))(|u_0|^p+|u_0|^q)\,dx +o(1) \\
&= \frac{1}{p}\int_{\R^N}(I_\alpha*|w_j|^p)|w_j|^p\,dx +\frac{1}{q}\int_{\R^N}(I_\alpha*|w_j|^q)|w_j|^q \,dx +o(1).
\end{aligned}
\end{equation}
Combining inequalities \eqref{HLS1} and \eqref{HLS2} with this, 
\[
\begin{aligned}
&\mathcal{S}_1\left(\int_{\R^N}(I_\alpha*|w_j|^p)|w_j|^p\right)^{\frac1p}
+\mathcal{S}_2\left(\int_{\R^N}(I_\alpha*|w_j|^q)|w_j|^q\right)^{\frac1q} \\
&\qquad\qquad\leq  \frac{1}{p}\int_{\R^N}(I_\alpha*|w_j|^p)|w_j|^p\,dx +\frac{1}{q}\int_{\R^N}(I_\alpha*|w_j|^q)|w_j|^q \,dx +o(1).
\end{aligned}
\]
We define
\[
x\coloneqq  \limsup_{j\to\infty} \int_{\R^N}(I_\alpha*|w_j|^p)|w_j|^p, \quad y \coloneqq \limsup_{j\to\infty} \int_{\R^N}(I_\alpha*|w_j|^q)|w_j|^q,
\]
both of which are finite since $\|w_j\|_{H^1}$ is bounded. Passing to a limit, we have
\begin{equation}\label{equ1}
\mathcal{S}_1x^{\frac1p}+\mathcal{S}_2y^{\frac1q} \leq \frac{1}{p}x+\frac{1}{q}y.
\end{equation}
We claim that $x = y = 0$. We prove this by getting rid of any other possibilities: $(1)\, x = 0, y \neq 0$;\,$(2)\, x\neq 0, y = 0$; and $(3)\, x\neq 0, y \neq 0$.
Suppose first the case $(1)$. Then one has 
$\mathcal{S}_2y^{\frac1q} \leq \frac{1}{q}y$,
which implies $y \geq (q\mathcal{S}_2)^{\frac{N+\alpha}{2+\alpha}}$.
In the case (2), we have $x \geq (p\mathcal{S}_1)^{\frac{N+\alpha}{\alpha}}$.
In the case (3), one has either  $y \geq (q\mathcal{S}_2)^{\frac{N+\alpha}{2+\alpha}}$ or  $x \geq (p\mathcal{S}_1)^{\frac{N+\alpha}{\alpha}}$
because, if $y < (q\mathcal{S}_2)^{\frac{N+\alpha}{2+\alpha}}$ and  $x < (p\mathcal{S}_1)^{\frac{N+\alpha}{\alpha}}$,
then 
\[
\frac1px+\frac1qy-\mathcal{S}_1x^{\frac1p}-\mathcal{S}_2y^{\frac1q} 
=   x^{\frac1p}(\frac1px^{1-\frac1p}-\mathcal{S}_1)+y^{\frac1q}(\frac1qy^{1-\frac1q}-\mathcal{S}_2) < 0
\]
so that \eqref{equ1} does not hold. 
Thus we conclude that in any case, either  $y \geq (q\mathcal{S}_2)^{\frac{N+\alpha}{2+\alpha}}$ or  $x \geq (p\mathcal{S}_1)^{\frac{N+\alpha}{\alpha}}$.
Now we again use Proposition \ref{Brezis-Lieb}, Proposition \ref{convergence} and \eqref{equ2} to deduce
\[
\begin{aligned}
J_\alpha(u_j) &= J_\alpha(u_0) +\frac{1}{2}\|w_j\|_{H^1}^2-\frac{1}{2p^2}\int_{\R^N}(I_\alpha*|w_j|^p)|w_j|^p\,dx \\
&\qquad\qquad -\frac{1}{2q^2}\int_{\R^N}(I_\alpha*|w_j|^q)|w_j|^q \,dx +o(1) \\
&=J_\alpha(u_0) +(\frac{1}{2p}-\frac{1}{2p^2})\int_{\R^N}(I_\alpha*|w_j|^p)|w_j|^p\,dx \\
&\qquad\qquad +(\frac{1}{2q}-\frac{1}{2q^2})\int_{\R^N}(I_\alpha*|w_j|^q)|w_j|^q \,dx +o(1). 
\end{aligned}
\]
Taking a limit and using the fact $J_\alpha(u_0) \geq 0$, we conclude that either
$c \geq (\frac{1}{2p}-\frac{1}{2p^2})(p\mathcal{S}_1)^{\frac{N+\alpha}{\alpha}}$ or
$c \geq (\frac{1}{2q}-\frac{1}{2q^2})(q\mathcal{S}_2)^{\frac{N+\alpha}{\alpha+2}}$ but this contradicts with the assumption of the proposition.
So the claim is proved. 

Now we are ready to complete the proof. Since $x = y = 0$, the equality \eqref{equ2} says that $\|w_j\|_{H^1} \to 0$ as $j \to \infty$ up to a subsequence. 
Therefore $u_j \to u_0$ in $H^1$ as $j \to \infty$ up to a subsequence.

\section {Proof of Theorem \ref{doubly-critical}}

We first show that $J_\alpha$ satisfies the mountain pass geometry on $H^1_r(\R^N)$. In other words, we show that there exist $r_0 > 0$ and $u_0 \in H^1_r(\R^N)$ such that
\begin{enumerate}[$(i)$]
\item $\displaystyle \inf_{\|u\|_{H^1_r} \leq r_0}J_\alpha(u) \geq 0 \quad \text{and} \quad \inf_{\|u\|_{H^1_r} = r_0}J_\alpha(u) > 0$, 
\item $\displaystyle J_\alpha(u_0) < 0 \quad (\text{and thus } \|u_0\|_{H^1_r} > r_0)$.
\end{enumerate}
The assertion $(i)$ immediately follows from the Hardy-Littlewood-Sobolev inequality (Proposition \ref{HLS}). 
Also, for any given $u \in H^1_r(\R^N)$, the function $f(t):= J_\alpha(tu)$ takes the form $At^2-Bt^{2p}-Ct^{2q}-Dt^{p+q}$ so that $f(0) = 0$ and $\lim_{t\to\infty}f(t) = -\infty$.
On the interval $(0,\,\infty)$, we can see that $f'(t) = 0$ if and only if $2pBt^{2p-2}+2qCt^{2q-2}+(p+q)Dt^{p+q-2} = 2A$.
Define $g(t):=2pBt^{2p-2}+2qCt^{2q-2}+(p+q)Dt^{p+q-2}$. Observe $g$ is strictly increasing on $(0,\,\infty)$, $g(0) = 0$ and $\lim_{t\to\infty}g(t) = \infty$.  
This shows $f$ admits a unique critical point $t_0$ on $(0,\,\infty)$ such that $f$ takes the maximum at $t = t_0$, $f$ is strictly increasing on $(0,\,t_0)$ and $f$ is strictly decreasing on $(t_0,\, \infty)$
and the assertion $(ii)$ follows.

Let $\Gamma$ denote the set of every continuous paths $\gamma:[0,\,1] \to H^1_r(\R^N)$ satisfying $\gamma(0) = 0$ and $J_\alpha(\gamma(1)) < 0$.    We define 
\[
c_0 = \inf_{u \in \Gamma}\max_{t \in [0,\,1]} J(\gamma(t)).
\] 
Since $J_\alpha$ enjoys the mountain pass geometry, the standard deformation lemma (see \cite{St,W}) says that there exists a (PS) sequence $\{u_j\} \subset H^1_r(\R^N)$ 
of $J_\alpha|_{H^1_r(\R^N)}$ at level $c_0$, i.e.,
\[
J_\alpha'(u_j) \to 0 \text{ in } H^{-1}_r(\R^N) \quad \text{and} \quad J_\alpha(u_j) \to c_0 \quad \text{as } j \to \infty.
\]
We claim that
\[
c_0 < \min\left(\frac{1}{2}(1-\frac{1}{p})p^{\frac{1}{p-1}}\mathcal{S}_1^{\frac{p}{p-1}},\, \frac{1}{2}(1-\frac{1}{q})q^{\frac{1}{q-1}}\mathcal{S}_2^{\frac{q}{q-1}}\right).
\]
If this is shown, it follows that $c_0$ is a critical level of $J_\alpha|_{H^1_r(\R^N)}$ by Proposition \ref{compact-level}. 
Since the mountain pass geometry of $J_\alpha$ implies $c_0 > 0$, we get a nontrivial critical point of $J_\alpha|_{H^1_r(\R^N)}$.
By the principle of symmetric criticality by Palais \cite{Pa}, this is also a nontrivial critical point of $J_\alpha$ on $H^1(\R^N)$, which is a solution to \eqref{main-eq}. 

Now, we show the claim holds. 
Let us define
\[
\mu_\lambda(x) \coloneqq A\frac{\lambda^{N/2}}{(\lambda^2+|x|^2)^{N/2}}, \quad \nu_\lambda(x) \coloneqq B\frac{\lambda^{(N-2)/2}}{(\lambda^2+|x|^2)^{(N-2)/2}},
\]
which are the extremal functions of  the inequalities \eqref{HLS1} and \eqref{HLS2} respectively.
The constants $A$ and $B$ are chosen to satisfy 
\[
\int_{\R^N}|\mu_1|^2\,dx = \int_{\R^N}(I_\alpha*|\mu_1|^p)|\mu_1|^p\,dx, \quad 
\int_{\R^N}|\nabla\nu_1|^2\,dx = \int_{\R^N}(I_\alpha*|\nu_1|^q)|\nu_1|^q\,dx.
\]
Since $N \geq 5$, one has $\mu_\lambda, \nu_\lambda \in H^1_r(\R^N)$.
Let $t_\lambda > 0$ and $s_\lambda > 0$ be two values satisfying
\[
J_\alpha(t_\lambda\mu_\lambda) = \max_{t> 0}J_\alpha(t\mu_\lambda), \quad
J_\alpha(s_\lambda\nu_\lambda) = \max_{t> 0}J_\alpha(t\nu_\lambda).
\]
Let $\tilde{t}_\lambda$ and $\tilde{s}_\lambda$ be the numbers that satisfy $J_\alpha(\tilde{t}_\lambda\mu_\lambda) < 0$ and $J_\alpha(\tilde{s}_\lambda\nu_\lambda) < 0$.
We have seen $\tilde{t}_\lambda > t_\lambda$ and $\tilde{s}_\lambda > s_\lambda$ should hold. Then by defining $\gamma_1(t) := t\tilde{t}_\lambda\mu_\lambda$ and $\gamma_2(t) := t\tilde{s}_\lambda\nu_\lambda$,
we see that
\[
c_0 \leq \min\{~ \max_{t \in [0,\,1]}J_\alpha(\gamma_1(t)),\, \max_{t\in [0,\,1]}J_\alpha(\gamma_2(t))~\}
= \min\{~J_\alpha(t_\lambda\mu_\lambda),\, J_\alpha(s_\lambda\nu_\lambda)~\}. 
\]

We compute
\begin{equation}\label{equ3}
\begin{aligned}
0 =& \frac{d}{dt}\bigg|_{t = t_\lambda}J_\alpha(t\mu_\lambda) \\
=& \, t_\lambda\int_{\R^N}|\nabla \mu_\lambda|^2\,dx+t_\lambda\int_{\R^N}|\mu_\lambda|^2\,dx \\
& -\frac{t_\lambda^{2p-1}}{p}\int_{\R^N}(I_\alpha*|\mu_\lambda|^p)|\mu_\lambda|^p\,dx
-\frac{t_\lambda^{2q-1}}{q}\int_{\R^N}(I_\alpha*|\mu_\lambda|^q)|\mu_\lambda|^q\,dx \\
&-\frac{(p+q)t_\lambda^{p+q-1}}{pq}\int_{\R^N}(I_\alpha*|\mu_\lambda|^p)|\mu_\lambda|^q\,dx \\
=& \, t_\lambda\lambda^{-2}\int_{\R^N}|\nabla \mu_1|^2\,dx+t_\lambda\int_{\R^N}|\mu_1|^2\,dx \\
& -\frac{t_\lambda^{2p-1}}{p}\int_{\R^N}(I_\alpha*|\mu_1|^p)|\mu_1|^p\,dx
-\frac{t_\lambda^{2q-1}\lambda^{-Nq+N+\alpha}}{q}\int_{\R^N}(I_\alpha*|\mu_1|^q)|\mu_1|^q\,dx \\
&-\frac{(p+q)t_\lambda^{p+q-1}\lambda^{-\frac{N}{2}(p+q)+N+\alpha}}{pq}\int_{\R^N}(I_\alpha*|\mu_1|^p)|\mu_1|^q\,dx.
\end{aligned}
\end{equation}
Let $t_\infty \coloneqq \limsup_{\lambda\to\infty}t_\lambda$. Suppose that $t_\infty = \infty$.
Then dividing the both side of \eqref{equ3} by $t_\lambda$ and taking a limit $\lambda\to\infty$, we get a contradiction and thus $t_\infty < \infty$.
We again pass to a limit $\lambda\to\infty$ in \eqref{equ3} to obtain
\[
0 = t_\infty\int_{\R^N}|\mu_1|^2\,dx -\frac{t_\infty^{2p-1}}{p}\int_{\R^N}(I_\alpha*|\mu_1|^p)|\mu_1|^p\,dx 
= \left(t_\infty-\frac{t_\infty^{2p-1}}{p}\right)\int_{\R^N}|\mu_1|^2\,dx,
\]
which implies $t_\infty = p^{1/(2p-2)}$. 

Now, observe
\[
\begin{aligned}
J_\alpha(t_\lambda\mu_\lambda)
=& \, \frac{t_\lambda^2\lambda^{-2}}{2}\int_{\R^N}|\nabla \mu_1|^2\,dx+\frac{t_\lambda^2}{2}\int_{\R^N}|\mu_1|^2\,dx \\
& -\frac{t_\lambda^{2p}}{2p^2}\int_{\R^N}(I_\alpha*|\mu_1|^p)|\mu_1|^p\,dx
-\frac{t_\lambda^{2q}\lambda^{-Nq+N+\alpha}}{2q^2}\int_{\R^N}(I_\alpha*|\mu_1|^q)|\mu_1|^q\,dx \\
&-\frac{t_\lambda^{p+q}\lambda^{-\frac{N}{2}(p+q)+N+\alpha}}{pq}\int_{\R^N}(I_\alpha*|\mu_1|^p)|\mu_1|^q\,dx \\
\leq& \, \left(\frac{t_\lambda^2}{2}-\frac{t_\lambda^{2p}}{2p^2}\right)\int_{\R^N}|\mu_1|^2\,dx
+\frac{t_\lambda^2\lambda^{-2}}{2}\int_{\R^N}|\nabla \mu_1|^2\,dx \\
& -\frac{t_\lambda^{p+q}\lambda^{-\frac{N}{2}(p+q)+N+\alpha}}{pq}\int_{\R^N}(I_\alpha*|\mu_1|^p)|\mu_1|^q\,dx \\
\end{aligned}
\]
Note that the curve $f(t) \coloneqq (\frac{t^2}{2}-\frac{t^{2p}}{2p^2})\int_{\R^N}|\mu_1|^2\,dx$ attains its maximum at $t = t_\infty$.
This shows
\[
\left(\frac{t_\lambda^2}{2}-\frac{t_\lambda^{2p}}{2p^2}\right)\int_{\R^N}|\mu_1|^2\,dx 
\leq \frac{1}{2}(1-\frac{1}{p})p^{\frac{1}{p-1}}\mathcal{S}_1^{\frac{p}{p-1}}
\]
Since $\frac{N}{2}(p+q)-(N+\alpha) < 2$ if and only if $4+\alpha < N$, we deduce that for sufficiently large $\lambda > 0$
\[
 J_\alpha(t_\lambda\mu_\lambda) < \frac{1}{2}(1-\frac{1}{p})p^{\frac{1}{p-1}}\mathcal{S}_1^{\frac{p}{p-1}}.
\]

Similarly we have
\begin{equation}\label{equ4}
\begin{aligned}
0 =& \frac{d}{dt}\bigg|_{t = s_\lambda}J_\lambda(t\nu_\lambda) \\
=& \, s_\lambda\int_{\R^N}|\nabla \nu_1|^2\,dx+s_\lambda\lambda^2\int_{\R^N}|\nu_1|^2\,dx \\
& -\frac{s_\lambda^{2p-1}\lambda^{-(N-2)p+N+\alpha}}{p}\int_{\R^N}(I_\alpha*|\nu_1|^p)|\nu_1|^p\,dx
-\frac{s_\lambda^{2q-1}}{q}\int_{\R^N}(I_\alpha*|\nu_1|^q)|\nu_1|^q\,dx \\
&-\frac{(p+q)s_\lambda^{p+q-1}\lambda^{-\frac{N-2}{2}(p+q)+N+\alpha}}{pq}\int_{\R^N}(I_\alpha*|\nu_1|^p)|\nu_1|^q\,dx,
\end{aligned}
\end{equation}
from which we conclude that $s_0 \coloneqq \limsup_{\lambda\to0}s_{\lambda} = q^{1/(2q-2)}$ by arguing similarly. 
Then we see that
\[
\begin{aligned}
 J_\alpha(s_\lambda\nu_\lambda) \leq& \left(\frac{s_\lambda^2}{2}-\frac{s_\lambda^{2q}}{2q^2}\right)\int_{\R^N}|\nabla \nu_1|^2\,dx
+\frac{s_\lambda^2\lambda^2}{2}\int_{\R^N}|\nu_1|^2\,dx \\
&-\frac{(p+q)s_\lambda^{p+q-1}\lambda^{-\frac{N-2}{2}(p+q)+N+\alpha}}{pq}\int_{\R^N}(I_\alpha*|\nu_1|^p)|\nu_1|^q\,dx \\
<& \frac{1}{2}(1-\frac{1}{q})q^{\frac{1}{q-1}}\mathcal{S}_2^{\frac{q}{q-1}}
\end{aligned}
\]
for sufficiently small $\lambda$ since $-\frac{N-2}{2}(p+q)+N+\alpha = \frac{N+\alpha}{N}< 2$.
This completes the proof.

\bigskip

{\bf Acknowledgements.}
This work was supported by Kyonggi University Research Grant 2016.

\end{document}